\documentclass[a4paper,11pt]{amsart}

\usepackage[all]{xy}
\usepackage{amssymb,amsmath,eufrak,url}
\author[Florent Benaych-Georges]{Florent Benaych-Georges}\address{Florent Benaych-Georges, LPMA, UPMC Univ Paris 6, Case courier 188, 4, Place Jussieu, 75252 Paris Cedex 05, France. } \email{florent.benaych@gmail.com}
\title[Cycles of random permutations with restricted cycle lengths]{Cycles of random permutations with restricted cycle lengths}
\date{\today}

\newcommand{\Sy}{\mathfrak{S}}
\newcommand{\Po}{\operatorname{\mathbb{P}}}
\newcommand{\Pon}{\operatorname{\mathbb{P}}_n}

\newcommand{\Poiss}{\operatorname{Poiss}}

\newcommand{\ds}{\displaystyle}

\newcommand{\ninf}{\underset{n\to\infty}{\longrightarrow}}
\newcommand{\ssi}{if and only if }

\newcommand{\E}{\mathbb{E}}

\newcommand{\R}{\mathbb{R}}
\newcommand{\C}{\mathbb{C}}
\newcommand{\n}{\mathbb{N}}

\newcommand{\pro}{probability }

\newcommand{\f}{\frac}
\newcommand{\ff}{\frac{1}}
\newcommand{\lf}{\left}
\newcommand{\ri}{\right}

\newcommand{\st}{such that }
\newcommand{\la}{\lambda}

\newcommand{\ste}{\, ;\, }
\newcommand{\mc}{\mathcal }

\newcommand{\eps}{\varepsilon}

\newcommand{\bck}{\backslash}

\newtheorem{Th}{Theorem}[section]
\newtheorem{propo}[Th]{Proposition} 
 
\newtheorem{lem}[Th]{Lemma}

\newtheorem{rmq}[Th]{Remark}
\newtheorem{cor}[Th]{Corollary}

\newenvironment{pr}{\noindent {\bf Proof. }}{\ \ \ \hfill$\square$}
\newenvironment{prth}{\noindent {\bf Proof of the theorem. }}{\ \ \ \hfill$\square$}

\long\def\symbolfootnote[#1]#2{\begingroup
\def\thefootnote{\fnsymbol{footnote}}\footnote[#1]{#2}\endgroup}

\begin{document}
\maketitle

\symbolfootnote[0]{{\it MSC 2000 subject classifications.}  primary 20B30, 
60B15, 
secondary 60C05.} 
\symbolfootnote[0]{{\it Key words.} Random permutations, cycles of random permutations, random permutations with restricted cycle lengths, $A$-permutations.}

\begin{abstract} We prove some general results about the asymptotics of the distribution of the number of cycles of given length of a random permutation whose distribution is invariant under conjugation. These results were first established to be applied in a forthcoming paper \cite{benaych.words}, where we prove results about cycles of random permutations which can be written as free words in several independent random permutations. However, we also apply them here to prove asymptotic results about random permutations with restricted cycle lengths. More specifically,  for $A$ a set of positive integers, we consider a random permutation chosen uniformly among the permutations of   $\{1,\ldots, n\}$ which have all their cycle lengths in $A$, and then let $n$ tend to infinity. Improving slightly a recent result of Yakymiv \cite{yakymiv}, we prove that under a general hypothesis on  $A$,  the numbers of cycles with fixed lengths of this  random permutation are asymptotically independent and distributed according to  Poisson distributions. In the case where $A$ is finite, we prove that  the behavior of these random variables is completely different: cycles with length $\max A$ are predominant.
\end{abstract}


\section*{Introduction}  It is well known that if for all positive integers $n$, $\sigma_n$ is a random permutation  chosen uniformly among all permutations of $\{1,\ldots, n\}$ and if for all positive integers $l$, $N_l(\sigma_n)$ denotes the number of cycles of length $l$ in the decomposition of $\sigma_n$ as a product of cycles with disjoint supports, then for all $l\geq 1$, the joint distribution of the random vector $$(N_1(\sigma_n),\ldots, N_l(\sigma_n))$$ converges weakly, as $n$ goes to infinity, to $$\Poiss(1/1)\otimes\Poiss (1/2)\otimes \cdots \otimes \Poiss (1/l),$$where for all positive number $\la$, $\Poiss(\la)$ denotes the Poisson distribution with mean $\la$. 

The proof of this result is rather simple (see, e.g.  \cite{ds94,bat}) because the uniform distribution on the symmetric group is easy to handle. However, many other distributions on the symmetric group give rise to limit distributions for ``small cycles", i.e. for the number of cycles of given length. In the first section of this paper, we shall  prove a general theorem about the convergence of the distributions of the number of cycles of given length of random permutations   distributed according to measures which are invariant under conjugation (Theorem \ref{theo.EVG-1.28.11.06.PM}). This result plays a key role in a forthcoming paper \cite{benaych.words}, were we prove results about cycles of random permutations which can be written as free words in several independent random permutations with restricted cycle length. More precisely, in \cite{benaych.words},  Corollary 3.2 (thus also, indirectly,  Theorems 3.7 and 3.8) and Theorem 3.12  are consequences of   Theorem \ref{theo.EVG-1.28.11.06.PM}  or of Corollary \ref{theo.EVG-1} of the present paper.

In the second part of the paper, for $A$ set of positive integers,  we introduce $\Sy_n^{(A)}$ to be the set of permutations of $\{1,\ldots, n\}$ which have all their cycle lengths in $A$ (such permutations are sometimes called {\it $A$-permutations}). For all $n$ \st  $\Sy_n^{(A)}$ is nonempty, we consider a random permutation $\sigma_n$ chosen uniformly in $\Sy_n^{(A)}$.  

We first prove, as an application of our general result mentioned above,  that under certain hypothesis on an infinite set $A$,    the result presented in the first paragraph about uniform random permutations stays as true as it can (as long as we consider the fact that for all $l\notin A$, $N_l(\sigma_n)$ is almost surely null): for all $l\geq 1$, the distribution of the random vector \begin{equation}\label{raylamontagne.12.07}(N_k(\sigma_n))_{\substack{1\leq k\leq l, k\in A}}\end{equation}converges weakly, as $n$ goes to infinity in such a way that $\Sy_n^{(A)}$ is  non empty, to  \begin{equation}\label{raylamontagne.2.12.07}\ds\underset{\substack{1\leq k\leq l, k\in A}}{\otimes }\Poiss (1/k).\end{equation}Here, we shall mention that as the author published this work on arxiv, it was pointed out to him that proving this result under some slightly stronger hypothesis was exactly the purpose of  a very recent paper   \cite{yakymiv}. However, the method used in this article is different from the one we use here: it relies on an identity in law between the random vector of \eqref{raylamontagne.12.07} and a vector with law \eqref{raylamontagne.2.12.07} conditioned to belong to a certain set and on  some estimations provided by asymptotic behavior of generating functions. It is the approach of  analytic combinatorics, which provides a powerful machinery for the analysis of random combinatorial objects. The book  \cite{flajolet} offers synthetic presentation of these tools. It is possible that  the result presented in this paragraph can be deduced from chapter IX of this book, but our proof is very short, and the object of the present paper is overall to prove the general result presented above about random permutations whose distributions are invariant under conjugation.

Note that the result presented in the previous paragraph  implies that  the number of cycles of any given length ``stays finite" even though $n$ goes to infinity,  i.e. takes large values with a very small probability. Hence if $A$ is finite, such a result cannot be expected. We also study this case here, and prove that if one denotes $\max A$ by $d$, for all $l\in A$, $N_l(\sigma_n)/n^{l/d}$ converges in every $L^p$ space to $1/l$. As a consequence, the cycles with length $d$ will be predominant: the cardinality of the subset of $\{1,\ldots, n\}$ covered by the supports of cycles with length $d$ in such a random permutation is asymptotic to $n$, which means that the random permutation is not far away from having order $d$. This remark will appear to be very helpful in the study of words in independent such random permutations.  
\\
\\
{\bf Notation.} In this text, we shall denote by $\n$  the set of nonnegative integers,. For $n$ an integer, we shall denote $\{1,\ldots, n\}$ by $[n]$ and    the group of permutations of $[n]$ by $\Sy_n$. For $A$ set of positive integers, $\Sy_n^{(A)}$ denotes  the set of permutations of $[n]$,  all of whose cycles have length in $A$.   For $\sigma\in \Sy_n$ and $l$ a positive integer, we shall denote by $N_l(\sigma)$ the number of cycles of length $l$ in the decomposition of $\sigma$ as a product of cycles with disjoint supports. For $\la >0$, $\Poiss(\la)$ will denote the Poisson distribution with parameter $\la$. If $I$ is a set, $|I|$ shall denote its cardinality. \\
\\
{\bf Acknowledgments.} The author would like to thank Professor Vatutin for having mentioned to him the works of Yakymiv and Pavlov, and Professor Yakymiv for his support and some interesting remarks about a draft of this work. He also would like to thank an anonymous 
referee for many interesting suggestions.

\section{A general result about cycles of random permutations}\label{section.gen.result.rand.perm}
\subsection{Main results} The main results of this section   are Theorem \ref{theo.EVG-1.28.11.06.PM} and Corollary \ref{theo.EVG-1}. Both of them play a key role in the paper \cite{benaych.words}. Recall that for all $n$ integer,  a \pro measure $\Po$ on $\Sy_n$ is said to be {\it invariant under conjugation} if  for all $\sigma,\tau\in \Sy_n$, $\Po(\{\sigma\circ\tau\circ\sigma^{-1}\})=\Po(\{\tau\})$.

\begin{Th}\label{theo.EVG-1.28.11.06.PM}Let $\mc{N}$ be an infinite set of positive integers. Fix a positive integer $q$,  some positive integers $l_1<\cdots<l_q$  and some  \pro measures $\mu_1,\ldots, \mu_q$  on the set of positive integers. 
   Let, for  for each $n\in A$, $\Po_n$ be \pro measure on $\Sy_n$ which    is invariant under conjugation. Suppose that for all $p\geq 1$, for all $k=(k_1, \ldots, k_q)\in \n^q$ \st $k_1l_1+\cdots +k_ql_q=p$ and for all $\sigma\in \Sy_{p}$ which has $k_1$ cycles of length $l_1$, \ldots, $k_q$ cycles of length $l_q$, the sequence $$\f{n^{p}}{l_1^{k_1}\cdots l_q^{k_q}k_1!\cdots k_q!}\Po_n(\{\tau \in \Sy_n\ste\forall i=1,\ldots,p, \tau(i)=\sigma(i)\})$$ converges, 
 as $n\in \mc{N}$ tends to infinity, to a limit, denoted by $S_{k}$, 
 \st for all $r_1,\ldots, r_q\geq 0$, we have  \begin{equation}\label{sweetheartcome.28.11.06} \ds\sum_{k_1\geq  r_1}\cdots \sum_{k_q\geq  r_q}(-1)^{k_1-r_1+\cdots +k_q-r_q}{k_1\choose r_1}\cdots {k_q\choose r_q} S_{(k_1,\ldots,k_q)}= \ds\prod_{\substack{1\leq i\leq q}}\mu_i(r_i).\end{equation}
 Then, if, for all $n\in \mc{N}$, $\sigma_n$ is a random variable distributed according to $\Po_n$,  
  the law of $(N_{l_1}(\sigma_n),$\ldots, $N_{l_q}(\sigma_n))$ converges, as $n\in \mc{N}$ tends to infinity, to $\mu_1\otimes \cdots \otimes \mu_q.$
 \end{Th}

\begin{rmq}\label{limit.proba.event.28.1106.plaid.rmq} Note that the series of (\ref{sweetheartcome.28.11.06}) are not asked to converge absolutely.  We only ask the sequence $$\ds\sum_{\substack{k_1\geq r_1,\ldots, k_q\geq r_q\\ k_1-r_1+\cdots +k_q-r_q\leq n}}(-1)^{k_1-r_1+\cdots +k_q-r_q}{k_1\choose r_1}\cdots {k_q\choose r_q} S_{(k_1,\ldots, k_q)},$$ to tend to the right hand term of (\ref{sweetheartcome.28.11.06}) as $n$ tends to infinity.\end{rmq}

Theorem \ref{theo.EVG-1.28.11.06.PM} will be proved in Section \ref{14.1.09.1}. Let us now give its main corollary.

\begin{cor}\label{theo.EVG-1}Let $\mc{N}$ be an infinite set of positive integers. Let $A$ be a set of positive integers and let, for each $n\in \mc{N}$, $\sigma_n$ be a random element of $\Sy_n$, distributed according to a \pro measure which is invariant under conjugation.  Suppose that for all $p\geq 1$, for all $\sigma\in \Sy_p^{(A)}$, the probability of the event $$\{\forall m=1,\ldots, p, \sigma_n(m)=\sigma(m)\}$$ is asymptotic to $n^{-p}$ as $n\in \mc{N}$ goes to infinity. Then
  for any finite subset $K$ of $A$, the law of $(N_l(\sigma_n))_{\substack{l\in K}}$ converges, as $n\in \mc{N}$ goes to infinity, to $\underset{l\in K}{\otimes} \Poiss(1/l).$
 \end{cor}
 
 \begin{rmq}\label{Wolfgang+WuTang.11.06.2}Note that the reciprocal implication is false. 
Consider for example a random permutation $\sigma_n$ with law $\lf(1-\ff{n}\ri)\mc{U}+\ff{n}\delta_{Id}$, where $\mc{U}$ denotes the uniform law on $\Sy_n$. Then the probability of the event $\{\sigma_n(1)=1\}$ is asymptotic to   $2/n$ as $n$ tends to infinity.
\end{rmq}
 
\noindent{\bf Proof of Corollary \ref{theo.EVG-1}. }The proof is immediate, since clearly, if one fixes a finite family $ l_1<\cdots< l_q$ of elements of $A$, then Theorem \ref{theo.EVG-1.28.11.06.PM} can be applied with $\Po_n=\operatorname{Law}(\sigma_n)$ for all $n$, with  $\mu_1=\Poiss(1/l_1),$\ldots,  $\mu_q=\Poiss(1/l_q)$ and with the $S_k$'s given by $$\forall k_1,\ldots, k_q, \quad S_{(k_1,\ldots, k_q)}=\ff{l_1^{k_1}\cdots l_q^{k_q} k_1!\cdots k_q!}.$$\hfill $\square$

\subsection{Technical preliminaries to the proof of Theorem \ref{theo.EVG-1.28.11.06.PM}}
Before the proof of Theorem \ref{theo.EVG-1.28.11.06.PM}, we shall prove  Proposition \ref{probasdavoirDEScertainsnombres}, which is a kind of Bonferroni inequality for inclusion-exclusion.  The principle is not new, but we did not find this result in the literature. 

Let us first recall Theorem 1.8 of \cite{bol}.\begin{Th}\label{reduc.to.triv}Fix $n,N\geq 1$, $\la_1,\ldots, \la_n\in \R$, $I_1,J_1,\ldots, I_n,J_n$ subsets of $[N]$. Then in order to have $$\ds\sum_{k=1}^n\la_k\Po( (\cap_{i\in I_i}A_i)\cap (\cap_{i\in J_i}(\Omega\bck A_i)))\geq 0$$ for any family $A_1,\ldots, A_N$ of events in any \pro space $(\Omega,\Sigma, \Po)$, it suffices to prove it under the additional hypothesis that each of the  $A_i$'s is either $\emptyset$ or $\Omega$. 
\end{Th}

\begin{propo}\label{probasdavoirDEScertainsnombres}Consider a \pro space $(\Omega,\Sigma,  \Po)$, $q\geq 1$, finite sets $I_1,\ldots, I_q$  and, for all $i=1,\ldots, q$,  $(A_{i,j})_{j\in I_i}$ a finite family of events of $\Sigma$. Let us define the random vector $C=(C_1,\ldots, C_q)$ by, for $i=1,\ldots, q$ and $\omega\in \Omega$, $$C_i(\omega)=|\{j\in I_i\ste \omega\in A_{i,j}\}|.$$Let us also define, for $k=(k_1,\ldots, k_q)\in \n^q\backslash\{0\}$, 
$$S_k=\ds\sum_{\substack{J_1\subset I_1\\ |J_1|=k_1}}\cdots\sum_{\substack{J_q\subset I_q\\ |J_q|=k_q}} \Po(\cap_{l=1}^q\cap_{j\in J_l}A_{i,j})$$and $S_0=1$. Then for all $r=(r_1,\ldots, r_q)\in \n^q$, \begin{equation}\label{31.05.05.1.prime} \Po(C=r)=\ds\sum_{k_1=r_1}^{{|I_1|}}\cdots\sum_{k_q=r_q}^{|I_q|} (-1)^{k_1-r_1+\cdots+k_q-r_q}{k_1\choose r_1}\cdots{k_q\choose r_q} S_{(k_1,\ldots, k_q)}.\end{equation}Moreover, alternating inequalities of the following type are satisfied: for all $m\geq 0$ odd (resp. even), \begin{equation}\label{31.05.05.1} \Po(C=r)\geq\ds\sum
 (-1)^{k_1-r_1+\cdots+k_q-r_q}{k_1\choose r_1}\cdots{k_q\choose r_q} S_{(k_1,\ldots, k_q)}\quad\textrm{ (resp. $\leq$),}\end{equation}where the sum runs over all families $(k_1,\ldots, k_q)$ of nonnegative integers \st $r_1\leq k_1\leq |I_1|$, \ldots, $r_q\leq k_q\leq |I_q|$ and $k_1-r_1+\cdots+k_q-r_q\leq m$.
\end{propo}

\begin{pr} Firstly, note that the alternating inequalities, used for $m$ large enough, imply (\ref{31.05.05.1.prime}). So we are only going to prove the alternating inequalities.

One can  suppose that for each $i=1,\ldots, q$, $I_i=[n_i]$, with $n_i$ a positive integer.
As an application of the previous theorem, one can suppose every $A_{i,j}$ to be  either $\emptyset$ or $\Omega$.  In this case, for all $i=1,\ldots, q$, the random variable $C_i$ is constant, equal to the number $c_i$ of $j$'s \st $A_{i,j}=\Omega$, and for all $k=(k_1,\ldots, k_q)\in \n^q$, $$S_k=\ds
{c_1\choose k_1}\cdots {c_q\choose k_q}.$$

Hence for $(r_1,\ldots, r_q)=(c_1,\ldots, c_q)$, for all $m\geq 0$, $$\ds\sum_{\substack{k_1=r_1,\ldots,n_1\\ \vdots\\ k_q=r_q,\ldots,n_q\\ k_1-r_1+\cdots+k_q-r_q\leq m}} (-1)^{k_1-r_1+\cdots+k_q-r_q}{k_1\choose r_1}\cdots{k_q\choose r_q} S_{(k_1,\ldots, k_q)}$$ $$=\ds\sum_{\substack{k_1=c_1,\ldots,n_1\\ \vdots\\ k_q=c_q,\ldots,n_q\\ k_1-r_1+\cdots+k_q-r_q\leq m}} (-1)^{k_1-r_1+\cdots+k_q-r_q}{k_1\choose c_1}\cdots{k_q\choose c_q} {c_1\choose k_1}\cdots {c_q\choose k_q},$$ which is equal to $1$, i.e. to $ \Po(C=r).$

Now, consider  $(r_1,\ldots, r_q)\neq (c_1,\ldots, c_q)$. Then $ \Po(C=r)=0$ and we have to prove that the right-hand-side term in equation (\ref{31.05.05.1}) is either nonnegative or nonpositive according to whether $m$ is even or odd.
For all $m\geq 0$, $$\ds\sum_{\substack{k_1=r_1,\ldots,n_1\\ \vdots\\ k_q=r_q,\ldots,n_q\\ k_1-r_1+\cdots+k_q-r_q\leq m}} (-1)^{k_1-r_1+\cdots+k_q-r_q}{k_1\choose r_1}\cdots{k_q\choose r_q} S_{(k_1,\ldots, k_q)}$$ $$=\ds\sum_{\substack{k_1=r_1,\ldots,n_1\\ \vdots\\ k_q=r_q,\ldots,n_q\\ k_1-r_1+\cdots+k_q-r_q\leq m}} (-1)^{k_1-r_1+\cdots+k_q-r_q}{k_1\choose r_1}\cdots{k_q\choose r_q} {c_1\choose k_1}\cdots {c_q\choose k_q}$$ $$=\ds\sum_{\substack{k_1=r_1,\ldots,c_1\\ \vdots\\ k_q=r_q,\ldots,c_q\\ k_1-r_1+\cdots+k_q-r_q\leq m}} (-1)^{k_1-r_1+\cdots+k_q-r_q}{k_1\choose r_1}\cdots{k_q\choose r_q} {c_1\choose k_1}\cdots {c_q\choose k_q}.$$ If there exists $i$ \st $r_i>c_i$, then the previous sum is zero. In the other case, since for all $0\leq r\leq k\leq c$,  ${k\choose r}{c\choose k}={c\choose r}{c-r\choose l}$ for $l=k-r$, the previous sum is equal to $$\ds{c_1\choose r_1}\cdots {c_q\choose r_q}\sum_{\substack{l_1=0,\ldots,c_1-r_1\\ \vdots\\ l_q=0,\ldots,c_q-r_q\\ l_1+\cdots +l_q\leq m}} (-1)^{l_1+\cdots +l_q}{c_1-r_1\choose l_1}\cdots{c_q-r_q\choose l_q}.$$
So we have to prove that for all $d=(d_1,\ldots, d_q)\in \n^q\bck\{0\}$ and for all $m\in \n$, 
$$\ds Z(m,d):=(-1)^m\sum_{\substack{l_1=0,\ldots,d_1\\ \vdots\\ l_q=0,\ldots,d_q\\ l_1+\cdots +l_q\leq m}} (-1)^{l_1+\cdots +l_q}{d_1\choose l_1}\cdots{d_q\choose l_q}$$is nonnegative. Let us prove it by induction over  $d_1+\cdots +d_q\geq 1$.

If $d_1+\cdots +d_q=1$, then $$Z(m,d)=\begin{cases}1&\textrm{if $m=0$,}\\ 0&\textrm{if $m>0$,}\end{cases}$$ so the result holds.

Suppose the result to be proved to the rank $d_1+\cdots +d_q-1\geq 1$. First note that if $m=0$, then $Z(m,d)=1$, so the result holds. So let us suppose that $m\geq 1$. Since $d_1+\cdots +d_q\geq 2$, there exists $i_0$ \st $d_{i_0}\neq 0$. One can suppose that $i_0=q$. Using ${d_q\choose l_q}={d_q-1\choose l_q}+{d_q-1\choose l_q-1}$,  one has $$Z(m,d) =Z(m,(d_1,\ldots, d_{q-1}, d_q-1))+Z(m-1, (d_1,\ldots, d_{q-1}, d_q-1)),$$ which completes the proof of the induction, and of the proposition.\end{pr}


\subsection{Proof of Theorem \ref{theo.EVG-1.28.11.06.PM}}\label{14.1.09.1}
Before the beginning of the proof, let us introduce some notation. 
 Let $\mathfrak{C}_l(n)$ be the set of cycles of $[n]$ with length $l$. 
Let,  for all cycle $c$ of $[n]$, $$E_c(n)=\{\sigma\in \Sy_n\ste \textrm{$c$ appears in the cycle decomposition of $\sigma$}\}.$$

{\it Step I. } In order to prove the theorem, we fix  a family of nonnegative integers $(r_1,\ldots, r_q)$, and we prove that the probability of the event $$\lf\{\forall i = 1,\ldots, q, N_{l_i}(\sigma_n)=r_i\ri\}$$ converges, as $n$ goes to infinity,   to $$ \ds\prod_{\substack{1\leq i\leq q}} \mu_i(r_i) ,$$i.e. to \begin{equation}\label{limit.proba.event} \ds\sum_{k_1\geq  r_1}\cdots \sum_{k_q\geq  r_q}(-1)^{k_1-r_1+\cdots +k_q-r_q}{k_1\choose r_1}\cdots {k_q\choose r_q} S_{(k_1,\ldots,k_q)}.\end{equation}
With the notations introduced above,  we have to prove that 
\begin{equation}\label{event}\Po_n(\forall i=1,\ldots,q, \textrm{ exactly $r_i$ of the events of the family $(E_c(n))_{c\in \mathfrak{C}_{l_i}(n)}$ occur})\end{equation}converges, as $n$ goes to infinity,  to (\ref{limit.proba.event}). 

By (\ref{31.05.05.1.prime}), for all $n$, the probability  of (\ref{event}) is \begin{equation}\label{proba.event}\ds\sum_{k_1=r_1,\ldots, |\mathfrak{C}_{l_1}(n)|}\cdots \sum_{k_q=r_q,\ldots,|\mathfrak{C}_{l_q}(n)|}(-1)^{k_1-r_1+\cdots +k_q-r_q}{k_1\choose r_1}\cdots {k_q\choose r_q} S_{(k_1,\ldots,k_q)}(n),\end{equation} where we have defined $S_0(n)=1$ and for all $k=(k_1,\ldots, k_q)\in \n^q\bck\{0\}$,\begin{equation}\label{trucadem.1.3.6.06} S_k(n):=\ds\sum
 \Po_n(\underset{i\in [q]}{\cap}\underset{c\in J_i}{\cap}E_{c}(n)),
\end{equation}
the sum running over all families $(J_i)_{i\in [q]}$ \st for all $i$,  $J_i\subset \mathfrak{C}_{l_i}(n)$ and $ |J_i|=k_i$.

{\it Step II. }Let us fix $k=(k_1,\ldots, k_q)\in \n^q\bck\{0\}$ and compute  $\ds\lim_{n\to \infty}S_k(n)$.  Define $p=k_1  l_1+\cdots +k_q  l_q$ and consider  $\sigma\in S_p$ such that the decomposition in cycles of $\sigma$ contains $k_1$ cycles of length $l_1$, $k_2$ cycles of length $l_2$, \ldots, $k_q$ cycles  of length $l_q$. Then the invariance of $\Po_n$ by conjugation   allows us to claim that  $S_k(n)$ is equal to  $$\Po_n(\{\sigma\in\Sy_n\ste \forall i=1,\ldots, p, \sigma_n(i)=\sigma(i)\})$$ 
 times the number of sets $J$ of cycles of $[n]$  which consist exactly in $k_1$ cycles of length $l_1$, $k_2$ cycles of length $l_2$, \ldots, $k_q$ cycles of length $l_q$ \st  these cycles are pairwise disjoint.
Such a set $J$ is defined by a set of pairwise disjoint subsets of $[n]$,  which consists exactly of $k_1$ subsets of cardinality $l_1$, $k_2$ subsets of cardinality  $l_2$, \ldots, $k_q$  subsets of cardinality $l_q$, and for every of these subsets, by the choice of a cycle having the subset for support. Hence  there are exactly $$\underbrace{\f{n!}{(n-p)!l_1!^{k_1}l_2!^{k_2}\cdots l_q!^{k_q}}\ff{k_1!k_2!\cdots k_q!}}_{\textrm{counting the sets of pairwise disjoint subsets of $[n]$}}\underbrace{(l_1-1)!^{k_1} (l_2-1)!^{k_2}(l_3-1)!^{k_3}\cdots (l_q-1)!^{k_q}}_{\textrm{choice of the cycles}}$$ such sets $J$. 
So $$S_k(n)=\f{n!}{(n-p)!l_1^{k_1}l_2^{k_2}\cdots l_q^{k_q}}\ff{k_1!k_2!\cdots k_q!}\Po_n(\{\sigma\in\Sy_n\ste \forall i=1,\ldots, p, \sigma_n(i)=\sigma(i)\}).$$
Hence by hypothesis, $$\ds\lim_{n\to\infty}S_k(n)=S_{(k_1,\ldots, k_q)}.$$

{\it Step III. }Now, let us prove that the probability of the event of (\ref{event}) converges to  (\ref{limit.proba.event}). Fix $\eps>0$. Choose $m_0\geq 0$ \st for all $m\geq m_0$, the absolute value of $$\ds\sum_{\substack{k_1\geq r_1\\ \vdots\\ k_q\geq r_q\\ k_1-r_1+\cdots +k_q-r_q>m}}(-1)^{k_1-r_1+\cdots +k_q-r_q}{k_1\choose r_1}\cdots {k_q\choose r_q} S_{(k_1,\ldots, k_q)},$$is less than $\eps/2$.

By (\ref{31.05.05.1}), for all $m,m'\geq m_0$ \st $m$ is odd and $m'$ is even, the \pro of the event of (\ref{event}) is bounded from below by $$\ds\sum_{\substack{k_1=r_1,\ldots, |C_1(n)|\\ \vdots \\ k_q=r_q,\ldots,|C_q(n)|\\ k_1-r_1+\cdots +k_q-r_q\leq m}}(-1)^{r_1+k_1+\cdots +r_q+k_q}{k_1\choose r_1}\cdots {k_q\choose r_q} S_{(k_1,\ldots,k_q)}(n)$$and bounded from above by  $$\ds\sum_{\substack{k_1=r_1,\ldots, |C_1(n)|\\ \vdots \\ k_q=r_q,\ldots,|C_q(n)|\\ k_1-r_1+\cdots +k_q-r_q\leq m'}}(-1)^{k_1-r_1+\cdots +k_q-r_q}{k_1\choose r_1}\cdots {k_q\choose r_q} S_{(k_1,\ldots,k_q)}(n).$$ Hence for $n$ large enough, the \pro of the event of (\ref{event}) is bounded from below by $$\ds-\eps/2+\sum_{\substack{k_1\geq r_1\\ \vdots \\ k_q\geq r_q\\ k_1-r_1+\cdots +k_q-r_q\leq m}}(-1)^{k_1-r_1+\cdots +k_q-r_q}{k_1\choose r_1}\cdots {k_q\choose r_q} S_{(k_1,\ldots,k_q)}$$and bounded from above by  $$\ds\eps/2+\sum_{\substack{k_1\geq r_1\\  \vdots \\ k_q\geq r_q\\ k_1-r_1+\cdots+ k_q-r_q\leq m'}}(-1)^{k_1-r_1+\cdots +k_q-r_q}{k_1\choose r_1}\cdots {k_q\choose r_q} S_{(k_1,\ldots,k_q)},$$hence is $\eps$-close to the sum of  (\ref{limit.proba.event}).
 It completes the proof of the theorem.
\ \ \ \hfill$\square$

\section{Cycles of random permutations with restricted cycle lengths}\label{section.S_n^A}
First of all, let us recall that  for $n$ large enough, $\Sy_n^{(A)}$ is non empty \ssi $n$ is divided by the greatest common divisor of $A$ (see Lemma 2.3 of \cite{neagu05} for example). 

\subsection{Case where $A$ is infinite}\label{chloe.please.be.the.one.I.expect.11.06.1}

The following proposition is the analog of the result stated in the beginning of the introduction, in the case where the random permutation we consider is not anymore distributed uniformly on the symmetric group but on the set of permutations all of whose  cycles lengths fall  in $A$ (note that  in this case,  for all $k\notin A$, $N_k(\sigma_n)$ is almost surely null).  

\begin{propo}\label{propo.limite.SyA}Suppose that $A$ is a set of positive integers such that, if one denotes by $q$ the greatest common divisor of $A$ and by $u_n$ the quotient $\lf|\Sy_{qn}^{(A)}\ri|/(qn)!$, one has
\begin{equation}\label{sara.veut.kon.parte.17.12.07}\f{u_n}{u_{n-1}}\ninf 1.
\end{equation}
 We consider, for $n$ large enough,   a random permutation $\sigma_n$  which has the uniform distribution on $\Sy_{qn}^{(A)}$. Then for all $l\geq 1$, the distribution of the random vector $$(N_k(\sigma_n))_{\substack{1\leq k\leq l, k\in A}}$$ converges weakly, as $n$ goes to infinity, to  $$\ds\underset{\substack{1\leq k\leq l, k\in A}}{\otimes }\Poiss (1/k).$$ \end{propo}

Note also that this result implies that for all $l$, even for large values of $n$, every $N_l(\sigma_n)$ takes large values with a very small probability.

\begin{pr} By Corollary \ref{theo.EVG-1}, it suffices to prove that for all $p\geq 1$, for all $\sigma\in \Sy_p^{(A)}$,   the probability  of the event $\{\forall m=1,\ldots, p, \sigma_n(m)=\sigma(m)\}$ is asymptotic to $n^{-p}$ as $n$ goes to infinity. 
This probability is equal to 
$$\f{\lf|\{s\in \Sy_{qn}^{(A)}\ste \forall m=1,\ldots, p, s(m)=\sigma(m)\}\ri|}{\lf| \Sy_{qn}^{(A)}\ri|}=\f{\lf| \Sy_{qn-p}^{(A)}\ri|}{\lf| \Sy_{qn}^{(A)}\ri|},$$
hence the proposition follows from  \eqref{sara.veut.kon.parte.17.12.07}.
\end{pr}

\begin{rmq}1. This result improves Theorem 1 of \cite{yakymiv}, which states the same result under the slightly stronger hypothesis that $(u_n)$ is a sequence with regular variation with exponent in $(-1,0]$. However, the author did not find any   example where the hypothesis of this result are satisfied but the hypothesis of Theorem 1 of \cite{yakymiv} are not. 
\\
2. A number of examples of  classes of sets $A$ satisfying the hypothesis of this proposition holds can be fund in the list of examples following Theorem 2 of \cite{yakymiv1}. It holds for example if  $\lf|A\cap [n]\ri|/n\ninf 1$. More details are given in Theorem 3.3.1 of the book \cite{yakymiv.book}.
\end{rmq}

\subsection{Case where $A$ is finite}\label{chloe.please.be.the.one.I.expect.11.06.2}
We are going to prove the following result: 

\begin{Th} Suppose that $A$ is a finite set of positive integers, and denote its maximum by $d$. We consider, for all $n$ such that  $\Sy_n^{(A)}$ is non empty,  a random permutation $\sigma_n$  which has the uniform distribution on $\Sy_n^{(A)}$. Then for all $l\in A$, as $n$ goes to infinity in such a way that $\Sy_n^{(A)}$ is non empty, $\f{N_l(\sigma_n)}{n^{l/d}}$ converges in all $L^p$ spaces ($p\in [1,+\infty)$) to $1/l$.
\end{Th}

\begin{rmq} It would be interesting to  know if we have a dilation of the random variables of $N_l(\sigma_n)/n^{l/d}-1/l$  which has a non degenerate weak limit as $n$ goes to infinity. It seems possible that analytic combinatorics, as presented in  \cite{flajolet}, could provide a way to answer this question. \end{rmq}

 To prove this theorem, we shall need the following lemmas. Lemma \ref{14.11.06.1} is well known (see, for instance, Theorem 3.53 of \cite{bona}). Lemma \ref{14.11.06.2} is Lemma 3.6 of \cite{neagu05}. 

\begin{lem}\label{14.11.06.1}Let $p$ be the greatest common divisor of $A$. Then for all complex number $z$, one has $$\ds\sum_{n\geq 0}\f{\lf|\Sy_{pn}^{(A)}\ri|}{(pn)!}z^{pn}=\exp\lf(\sum_{k\in A}\f{z^k}{k}\ri).$$
\end{lem}

\begin{lem}\label{14.11.06.2} Let $B$ be a finite set of positive integers. 
Let $(c_j)_{j\in B}$ be a finite family of positive numbers. Let $\ds\sum_{n\geq1}b_nw^n$ be the power expansion   of $\exp \lf(\sum_{j\in B}c_jw^j\ri).$ Suppose that $b_n>0$ for sufficiently large $n$. Then, as $n$ goes to infinity, $$\f{b_{n-1}}{b_n}\sim\lf(\f{n}{bc_b}\ri)^{1/b},$$with $b=\max B$.
\end{lem}

\begin{prth} First note that by H\"older formula, it suffices to prove that for all $p$ positive integer, the expectation of the $2p$-th power of $$\f{N_l(\sigma_n)}{n^{l/d}}-\ff{l}$$ tends to zero as $n$ goes to infinity. Hence by the binomial identity, it suffices to prove that for all $l\in A$, for all $m\geq 1$, the expectation of the $m$-th power of $N_l(\sigma_n)$ is asymptotic to $n^{ml/d}/l^m$ as $n$ goes to infinity in such a way that $\Sy_n^{(A)}$ is non empty. 

One can suppose that for all such $n$, the \pro space where $\sigma_n$ is defined is $\Sy_n^{(A)}$, endowed with the uniform \pro measure $\Pon$. Let  $\E_n$ denote the expectation with respect to $\Pon$. 

So let us fix $l\in A$ and $m\geq 1$. Since $N_l(\sigma_n)=\ds\ff{l}\sum_{k=1}^n1_{\{\textrm{$k$ belongs to a cycle of length $l$}\}},$ one has
 $$\E_n[N_l(\sigma_n)^m]=\ds\ff{l^m}\sum_{\substack{m_1,\ldots, m_n\geq 0\\ m_1+\cdots +m_n=m}}{m\choose m_1, \ldots, m_n} \E_n\lf[
\prod_{k=1}^n(1_{\{\textrm{$k$ belongs to a cycle of length $l$}\}})^{m_k}
\ri]
$$
But  $\Pon$ is invariant by conjugation, so for all $m_1,\ldots, m_n\geq 0$, 
$$\E_n\lf[
\prod_{k=1}^n(1_{\{\textrm{$k$ belongs to a cycle of length $l$}\}})^{m_k}
\ri]
 $$
depends only on the number $j$ of $k$'s \st $m_k\neq 0$. So \begin{eqnarray}
\E_n[N_l(\sigma_n)^m]&=&\ds\ff{l^m}\sum_{j=1}^m\sum_{\substack{m_1,\ldots, m_n\geq 0\\ |\{k\in [n]\ste m_k\neq 0\}|=j\\ m_1+\cdots +m_n=m}}{m\choose m_1, \ldots, m_n} \Pon\lf(\textrm{$1,\ldots, j$ belong to cycles of length $l$}\ri)\nonumber\\
&=& \ds\ff{l^m}\sum_{j=1}^m\lf[{n \choose j} \Pon\lf(\textrm{$1,\ldots, j$ belong to cycles of length $l$}\ri) \!\!\!\!\!\!\sum_{\substack{m_1,\ldots, m_j\geq 1\\ m_1+\ldots +m_j=m}}{m\choose m_1, \ldots, m_j} \ri].\label{14.11.06.999}\end{eqnarray}

Now, let us fix   $j\geq 1$ and let us denote   by $\operatorname{P}(j)$ the set of partitions of $[j]$. We have
 \begin{eqnarray}&&\Pon\lf(\textrm{$1,\ldots, j$ belong to cycles of length $l$}\ri)\nonumber \\ &=&\sum_{\pi\in \operatorname{P}(j)}\Pon(\textrm{$1,\ldots, j$ are in cycles of length $l$}\nonumber\\ &&\quad\quad\quad \quad\quad{\textrm{and $\forall i,i'\in [j]$, [$i,i'$ belong to the same cycle] $\Leftrightarrow [i=i'\mod \pi]$})}\nonumber
\\ &=& \sum_{ \substack{ \pi\in \operatorname{P}(j)\\ \pi=\{V_1,\ldots, V_{|\pi|}\} } } 
{n-j \choose l-|V_1|,\ldots ,l-|V_{|\pi|}|, n-l|\pi|} \lf((l-1)!\ri)^{|\pi|} \f{ \lf|\Sy_{n-l|\pi|}^{(A)}\ri|}{  \lf|\Sy_{n}^{(A)}\ri|}\nonumber\\ 
&=& \sum_{\pi\in \operatorname{P}(j)} \ff{n(n-1)\cdots (n-j+1)}\f{ \lf|\Sy_{n-l|\pi|}^{(A)}\ri|/(n-l|\pi|)!}{  \lf|\Sy_{n}^{(A)}\ri|/n!}\prod_{V\in \pi}\f{(l-1)!}{(l-|V|)!}.\label{melon.beb's.s'en.va}\end{eqnarray}

Let $p$ be the greatest common divisor of $A$. We know  \cite[Lem. 2.3]{neagu05} that for all positive integer $n$, $\Sy_n^{(A)}\neq \emptyset\Longrightarrow p|n,$ and that for sufficiently large $n$, the inverse implication is also true. Hence by Lemma \ref{14.11.06.1}, for $z\in \C$, one has $$\ds\sum_{n\geq 0}\f{\lf|\Sy_{pn}^{(A)}\ri|}{(pn)!}(z^{p})^n=\exp\lf(\sum_{j\in \ff{p}\cdot A}\f{(z^p)^j}{pj}\ri).$$ Hence for $w\in \C$, one has $$\ds\sum_{n\geq 0}\f{\lf|\Sy_{pn}^{(A)}\ri|}{(pn)!}w^n=\exp\lf(\sum_{j\in \ff{p}\cdot A}\f{w^j}{pj}\ri).$$ So by Lemma \ref{14.11.06.2}, as $n$ goes to infinity, one has $$\f{\lf|\Sy_{pn-p}^{(A)}\ri|/(pn-p)!}{\lf|\Sy_{pn}^{(A)}\ri|/(pn)!}\sim\lf(\f{n}{(d/p)1/d}\ri)^{p/d}=\lf(pn\ri)^{p/d}$$(note that in the case where the greatest common divisor of $A$ is $1$,  this result can be deduced from the main theorem of \cite{pavlov95}).
It follows, by induction on $k$,   that, as $n$ goes to infinity in such a way that $p$ divides $n$,  for any positive integer $k$ divisible by $p$, we have
$$\f{\lf|\Sy_{n-k}^{(A)}\ri|/(n-k)!}{\lf|\Sy_{n}^{(A)}\ri|/(n)!} \sim n^{k/d}.$$
 Hence in \eqref{melon.beb's.s'en.va}, for each partition $\pi$, the term corresponding to $\pi$ is asymptotic to $$ n^{l|\pi|/d-j}\prod_{V\in \pi}\f{(l-1)!}{(l-|V|)!},$$ thus  in \eqref{melon.beb's.s'en.va},  the leading term is the one of the singletons partition,  and \begin{equation}\label{14.11.06.1000}
\Pon\lf(\textrm{$1,\ldots, j$ belong to cycles of length $l$}\ri)\sim n^{(l/d-1)j}.
\end{equation}

Combining (\ref{14.11.06.999}) and  (\ref{14.11.06.1000}), one gets $\E_n[N_l(\sigma_n)^m]\sim\ds\f{n^{lm/d}}{l^m}$, which completes the proof of the theorem.
\end{prth}

\end{document}